\documentclass[11pt,oneside]{article}
\usepackage{latexsym,amsmath,amssymb,amsfonts,amsbsy,fancyhdr}
\usepackage{theorem}
\usepackage{amscd}
\usepackage{enumerate}
\usepackage{mathrsfs}
\usepackage{dsfont}
\usepackage[active]{srcltx}
\usepackage[colorlinks=false, pdfstartview=FitV, linkcolor=red, citecolor=green, urlcolor=red]{hyperref}
\usepackage{nicefrac}
\usepackage{marvosym}
\usepackage{longtable}
\usepackage{subfigure}
\usepackage{setspace}

\newcommand{\R}{\mathbb R}

\newcommand{\N}{\mathbb N}

\newcommand{\E}{\mathbb E}

\renewcommand{\L}{\mathscr L}

\newcommand{\dd}{{\rm{d}}}

\newcommand{\du}{\frac{\partial}{\partial u}}
\newcommand{\ds}{\frac{\partial}{\partial s}}

\newtheorem{theorem}{Theorem}

\newtheorem{lemma}[theorem]{Lemma}
\theorembodyfont{\rm}
\newtheorem{definition}[theorem]{Definition}
\newtheorem{example}[theorem]{Example}

\newtheorem{remark}[theorem]{Remark}

\newtheorem{proof}{Proof}

\newenvironment{pf}{\begin{proof}}{\hfill$\square$\end{proof}}

\allowdisplaybreaks

\DeclareMathAlphabet{\mathcal}{OMS}{cmsy}{m}{n}

\begin{document}
\title{Maximal inequalities for fractional L\'evy and related processes}
\author{Christian Bender, Robert Knobloch and Philip Oberacker
\thanks{Department of Mathematics, Saarland University, P.O. Box 151150, 66041
Saarbr\"ucken, Germany
\newline
e-mail: bender@math.uni-sb.de, knobloch@math.uni-sb.de,
oberacker@math.uni-sb.de}}
\maketitle

\begin{abstract}  
In this paper we study processes which are constructed by a convolution of a deterministic kernel with a martingale. A special emphasis is put on the case where the driving martingale is a centred 
L\'evy process, which covers the popular class of fractional L\'evy processes. As a main result we show that, under appropriate assumptions on the kernel and the martingale, the maximum process of the corresponding `convoluted martingale' is $p$-integrable and we derive maximal inequalities in terms of the kernel and of the moments of the driving martingale.   
\end{abstract}


\section{Introduction}
 
We study a class of stochastic processes which is defined via a convolution of a deterministic Volterra kernel with a stochastic process. Precisely, given a two-sided martingale $X$
and a continuous deterministic kernel $f$ we consider the process 
$$
M(t)=\int_{-\infty}^t f(t,s)X(ds),\;0\leq t \leq T.
$$
Here the behaviour of the kernel $f$ and the martingale $X$ at minus infinity must be balanced in an appropriate way to ensure that this integral exists (at least in an improper sense).
A particular emphasis will be given to the case where the martingale is a centred L\'evy process, in which case we refer to the above class of processes as L\'evy-driven Volterra processes. An intriguing feature of L\'evy-driven Volterra processes from the modelling point of view is that the second order structure (and, hence, the memory) is encoded in the kernel,
while other distributional properties such as the tail behaviour can be goverened by the choice of the driving process. L\'evy-driven Volterra processes and related models have thus been applied to various problems in mathematical finance,
see e.g. \cite{BBV_mesp, B_mqle, F_ccf, KM_gfl}, but they are also of interest in other fields such as network traffic \cite{WT05} or signal processing~\cite{Ual14}.

In this paper we are mainly concerned with the question under which conditions on the kernel integrability properties of the driving martingale $X$ are inherited by the convoluted process 
$M$. More precisely, we derive $\mathscr L^p(\mathbb P)$-inequalities for the maximum process $M^*(T):=\sup_{s\in [0,T]} |M(s)|$ in terms of the kernel and the 
$\mathscr L^p(\mathbb P)$-norm of the driving martingale. Let us illustrate our results for the case of a fractional L\'evy process
$$
M_d(t)=\int_{-\infty}^t  \frac{1}{\Gamma(d+1)}\left((t-s)^d_+-(-s)^d_+\right) L(ds),
$$
where $d\in (0,1/2)$ and $L$ is a centred two-sided L\'evy process with a finite second moment. We will show in Example \ref{ex_fLp} below that for every $p\geq 2$ and $\delta>0$ such that  $d+\delta<1/2$ there is a constant
  $C_{p,\delta,d}$ independent of the driving L\'evy process $L$ such that for every $T\geq 1$
$$
\E\left( \sup_{0\leq t \leq T} |M_d(t)|^p\right) \leq C_{p,\delta,d} \E(|L(1)|^p) \; T^{p(d+1/2+\delta)}.
$$
As $\E (|M_d(T)|^2)= \E (|M_d(1)|^2) T^{2(d+1/2)}$ by \cite{Mar06}, we observe that $p(d+1/2)$ is a lower bound for the optimal rate 
in $T$ in the above maximal inequality and our rate is arbitrarily close to this expected optimal rate.
 We also stress that in the case of a fractional Brownian motion the rate $p(d+1/2)$ is an obvious consequence of the self-similarity and refer to \cite{NV_mifbm} for further discussion on maximal
inequalities for fractional Brownian motion. However, square integrable fractional L\'evy processes are not self-similar except in the fractional Brownian motion case (see \cite{Mar06}) and hence 
this line of reasoning cannot be applied in the context of the present paper. 

The paper is organised as follows: In Section 2 we treat the case of a general martingale $X$ as driving process. We introduce a suitable class of kernels which ensures 
that the convoluted process exists and has a c\`adl\`ag modification. More precisely, we also relate the jumps of the convoluted process to the jumps of the driving martingale. We also derive a first version of an $\mathscr L^p(\mathbb P)$-maximal inequality which depends on the asymptotic behaviour of 
the function $\E(|X(2t)-X(t)|^p)$ as $t$ approaches infinity. In Section 3 the results are then refined for the case of a driving centred L\'evy process exploiting the stationary increments.
This leads to an improved maximal inequality for L\'evy-driven Volterra processes. We consider the situation of a L\'evy process with finite second moment in Theorem~\ref{t.1_levy}
and that of a symmetric $\alpha$-stable process in Theorem~\ref{t.2_levy}.

\section{The general case}\label{s.main_results}
In this section we derive paths properties and a maximal inequality for processes of the form

\[
 M(t) =  \int_{-\infty}^t f(t,s) \ X(\dd s),\; 0\leq t \leq T,
\]

where $M$ is a two-sided martingale, $f$ is a deterministic kernel function and $T>0$ is fixed.

\begin{definition}\label{d.X}
Let $\hat X := (\hat X(t))_{t \ge0}$ be a c\`adl\`ag martingale starting at zero. We construct a two-sided process  $X:=(X(t))_{t \in \R}$ by taking two independent copies $(X_1(t))_{t \geq 0}$ and $(X_2(t))_{t \geq 0}$ of  $\hat{X}$ and defining

\begin{equation}\label{e.two-sided}
 X(t) := \begin{cases} X_1(t), &  t \geq 0\\-X_2(-(t-)), & t<0.
 \end{cases}
\end{equation}
\end{definition}

Throughout the paper $\varphi:[0,\infty)\to [1,\infty)$ denotes a nondecreasing function. Let us now introduce the following class of Volterra type kernels depending on $\varphi$:
\begin{definition}\label{d.1}
Let  $\tau \in [-\infty, 0]$.
We denote by $\mathcal K(\varphi, \tau)$ the class of measurable functions $f:\R^2\to\R$ with $\text{supp}f \subset [\tau, \infty)^2$ such that
\begin{enumerate}[(i)]
%\item\label{support} there exists some $\tau \in [-\infty, 0]$ such that $\text{supp}f \subset [\tau, \infty)^2$,
\item\label{volterra} $\forall\,s>t\ge0:\quad f(t,s)=0$,

\item\label{continuity_of_f} %for every $s\in\R$ the function $f(\cdot,s)$ is %continuous c\`adl\`ag on $[s,T]$,
%for every $s\in\R$ 
the mapping $t\mapsto f(t,t)$ is continuous on $[0,T]$; moreover if $\tau>-\infty$ then also $t\mapsto f(t,\tau)$ is continuous  on $[0,T]$,

\item for all $t \in [0,T]$ we have
    
    \begin{equation}\label{e.decay_of_f}
     \lim_{s \to \infty} f(t,-s) \varphi(s) =0,
    \end{equation}

    \item\label{derivative_second_argument} for every fixed $t \in [0,T]$ the function $s \mapsto f(t,s)$ is absolutely continuous on $[\tau,t]$ with density $\ds f(t,\cdot)$, i.e.
    \[
     f(t,s) = f(t,\tau) +\int_\tau^s \du f(t,u) \ \dd u, \qquad \tau \le s \le t,
    \]
    where $f(t,-\infty) := \lim_{x \to -\infty} f(t,x) = 0$, such that
    \begin{enumerate}
    \item the function $t \mapsto \ds f(t,s)$ is continuous on $(s,\infty)$ for $\lambda$-a.e. $s \in [\tau,\infty)$, where $\lambda$ denotes the Lebesgue measure,
    \item there exists an $\epsilon >0$ (independent of $t$) such that
    
    \begin{equation}\label{e.sup_ds}
     \sup_{t \in [0,T]} \int_{-\infty}^t \left| \ds f(t,s) \varphi(|s|) \right|^{1+\epsilon} \left(|s|^{2\epsilon}  \lor 1 \right) \ \dd s < \infty.
    \end{equation}

    \end{enumerate}
    \end{enumerate}
\end{definition}  

The function $\varphi$ describes the behaviour of the kernel $f$ and its density at $s=-\infty$. If it is connected to the $\mathscr L^p(\mathbb P)$-norm, denoted by $\|\cdot\|_p$, of the increments of the martingale $X$ in 
an appropriate way, we will show that the improper integral in the definition of $M$ exists and that $M$ inherits path properties and finite moments from $X$.

\begin{theorem}\label{t.0}
Let $f\in\mathcal K(\varphi, \tau)$, $p>1$ with $X(t)\in\mathscr L^p(\mathbb P)$ for every $t \in \mathbb{R}$ and assume that

\begin{equation}\label{e.int_X_phi_finite}
\sum_{n=0}^\infty \frac{\left\| X( 2^{n+1})- X(2^{n})\right\|_p}{\varphi(2^{n})} \ < \infty.
\end{equation}

Then the following assertions hold:

\begin{enumerate}
 \item The limit

\begin{equation}\label{d.2}
\tilde M(t):=\lim_{n \to \infty} \int_{-n}^t f(t,s) \ X(\dd s),\quad 0\leq t \leq T,
\end{equation}

exists $\mathbb P$-a.s. and in $\mathscr L^p(\mathbb P)$ and a modification of $\tilde M$ is given by

\begin{equation}\label{integration_by_parts_2}
  M(t) := f(t,t) X(t) -f(t,\tau) X(\tau) - \int_{\tau}^t X(s) \ds f(t,s) \ \dd s,
\end{equation}

where $f(t,-\infty) X(-\infty) := \lim_{N \to \infty} f(t,-N) X(-N)=0$ holds. 

\item The process $M$ has c\`adl\`ag paths and $\Delta M(t)=f(t,t)\Delta X(t)$.

\item The following maximal inequality holds:
\begin{equation}\label{t.1}
 \begin{aligned}
& \left\| \sup_{t\in[0,T]}|M(t)| \,\right\|_p  
\\ 
\le &\; \frac{p}{p-1} \sup_{t\in[0,T]}|f(t,t)|  \|X(T)\|_p
+ \sup_{t\in[0,T]} \|X(\tau)f(t,\tau)\|_p 
\\ 
& +  \frac{2p}{p-1}  \sup_{t\in[0,T]} \left(\int_{\tau}^t \varphi(|s|) \left|\ds f(t,s)\right| \ \dd s\right) \left(
\|X(1)\|_p+\sum_{n=0}^\infty \frac{\left\| X( 2^{n+1})- X(2^{n})\right\|_p}{\varphi(2^{n})}  \right)
\end{aligned}
\end{equation}
where $f(t,-\infty)X(-\infty)=0$ (cf. 1.).
%\begin{equation}\label{t.1}
% \sup_{t \in [0,T]}|M(t)|\in\mathscr L^p(\mathbb P).
%\end{equation}
\end{enumerate}
\end{theorem}

\begin{remark}\label{r.integral}
 Note that by an application of H\"older's inequality we have
\begin{align*}
 &\sup_{t \in [0,T]}\int_{\tau}^t \varphi(|s|) \left|\ds f(t,s)\right| \ \dd s 
 \\
 &\leq \left(\sup_{t \in [0,T]}\int_{-\infty}^t \left| \ds f(t,s) \varphi(|s|) \right|^{1+\epsilon} \left(|s|^{2\epsilon}  \lor 1 \right) \ \dd s\right)^{\frac{1}{1+\epsilon}}
\left(\int_{-\infty}^T \left(|s|^{-2}  \wedge 1 \right) \ \dd s\right)^{\frac{\epsilon}{1+\epsilon}}
\end{align*} 
and hence the right-hand side of (\ref{t.1}) is finite for $f\in K(\varphi, \tau)$, because of \eqref{e.sup_ds}.
\end{remark}

\begin{remark}
 \begin{itemize}
  \item If $f$ has compact support (i.e. $\tau>-\infty$), one can always choose $\varphi\equiv 1$ and consider kernels $f$ in the class $K(1, \tau)$. 
Indeed, we can replace $X(t)$ by $X(-(T\wedge |\tau|)\vee t \wedge (T\wedge |\tau|))$,
which does not change the definition of $M$, but ensures that (\ref{e.int_X_phi_finite}) is satisfied.
\item If $f$ is sufficiently regular and has compact support, the above relation between the jumps of $M$ and $X$ has already been observed in several papers e.g.
\cite{BM_sccl} and \cite{MN_sppvp}. Without such regularity assumptions $M$ may fail to be continuous, even if $f$ vanishes on the diagonal. This has been shown 
by a counterexample in \cite{KMR_trcbsc}.
 \end{itemize}

\end{remark}

We now prepare the proof of Theorem \ref{t.0} by the following lemma.

\begin{lemma}\label{e.exp.sup.finite}
 Let $p>1$ be as in Theorem~\ref{t.0}. Then
 
  \[
   \left\|\sup_{s \in \R} \frac{|X(s)|}{\varphi(|s|)^{}} \right\|_p  \le \frac{2p}{p-1} \left( \|X(1)\|_p + \sum_{n=0}^\infty \frac{\left\| X(2^{n+1})- X(2^{n})\right\|_p}{\varphi(2^{n})}  \right)< \infty.
  \]

\end{lemma}

\begin{pf}
Since Lemma~\ref{e.exp.sup.finite} is only concerned with a distributional property of $X$, the construction of the two-sided process $X$ 
entails that it suffices to consider $X$ on the positive half line.

We first introduce the abbreviation
 
 \[
  \mathfrak L_N := \sup_{s \in [0,2^N]} \frac{|X(s)|}{\varphi(s)^{}} \quad (N \in \N).
 \]

 By using Doob's inequality and the fact that the mapping $s \mapsto \varphi(s)^{-1}$ is bounded by 1, we 
infer $\mathfrak L _N\in \L^p(\mathbb P)$. Drawing a distinction whether the supremum in the expression $\mathfrak L_N$ is attained on the set $[0,2^{N-1}]$ or   $[2^{N-1}, 2^N]$ and using the reverse triangle inequality as well as the  fact that $\varphi$ is nondecreasing  we continue with the chain of estimates
 
 %\begin{equation}\label{e.6}
 \[
  \begin{aligned}
   \mathfrak L_N &= (\mathfrak L_N- \mathfrak L_{N-1}) + \mathfrak L_{N-1} \le \left( \sup_{s \in [2^{N-1},2^N]} \frac{|X(s)|}{\varphi(s)^{}}- \frac{|X(2^{N-1})|}{\varphi(2^{N-1})^{}}\right)+  \mathfrak L_{N-1}\\
   &\le \sup_{s \in [2^{N-1},2^N]} \frac{\left| X(s) - X(2^{N-1})\right|}{\varphi(2^{N-1})^{}}+  \mathfrak L_{N-1}.
  \end{aligned}
 \]

 Proceeding inductively we obtain
 
 %\begin{equation}\label{e.7}
 \[
    \mathfrak L_N \le \sum_{n=1}^N \frac{\sup_{s \in [2^{n-1},2^n]} \left| X(s)- X(2^{n-1})\right|}{\varphi(2^{n-1})^{}}+ \sup_{s \in [0,1]} |X(s)|.
 \]

We now use Minkowski's inequality and Doob's inequality to deduce
 
  %\begin{equation}\label{e.8}
  \[
  \begin{aligned}
   \left\|\mathfrak L_N\right\|_p &\le  \left\|\sup_{s \in [0,1]} |X(s)| \right\|_p + \sum_{n=1}^N \frac{\left\|\sup_{ s \in [0,2^{n-1}]} |X(s+2^{n-1})- X(2^{n-1})|\,\right\|_p}{\varphi(2^{n-1})^{}} \\
   &\le   \frac{p}{p-1} \left(\left\|X(1)\right\|_p+\sum_{n=1}^N \frac{\left\| X(2 \cdot 2^{n-1})- X(2^{n-1})\right\|_p}{\varphi(2^{n-1})^{}} \right).
  \end{aligned}
 \]
 
 Since the series $\sum_{n=1}^\infty \frac{\left\| X(2^n)- X(2^{n-1})\right\|_p}{\varphi(2^{n-1})}$ is finite by assumption \eqref{e.int_X_phi_finite}, it thus follows from the monotone convergence theorem that
 \[
  \begin{aligned}
   \left\|\sup_{s \in [0, \infty)} \left( |X(s)| \varphi(s)^{-1} \right)\right\|_p  &= \left\|\lim_{N \to \infty} \mathfrak L_N\right\|_p = \lim_{N \to \infty}  \left\|\mathfrak L_N\right\|_p \\
   &\le  \frac{p}{p-1} \left(\left\|X(1)\right\|_p+ \sum_{n=1}^\infty \frac{\left\| X(2^n)- X(2^{n-1})\right\|_p}{\varphi(2^{n-1})} \right)\\
  &< \infty.
  \end{aligned}
\]

By symmetry the same inequality holds for the negative half line, which completes the proof.
\end{pf}

We are now in a position to prove Theorem~\ref{t.0}.

{\bf Proof of Theorem~\ref{t.0}}
%
%We prepare the proof with the following consideration.  
%In particular,
%\begin{eqnarray}\label{L1}
% \int_{\tau}^t \left| \ds f(t,s)\varphi(|s|)\right| \ \dd s = \int_{\tau}^t \left| \Upsilon(t,s) \right| \frac{1}{1\lor|s|^\gamma}  \ \dd s<\infty
%\end{eqnarray}
%
%

1. For every $n \in \N$ we define
 \[
  M_n(t) := \int_{\tau \lor -n}^t f(t,s) \ X( \dd s ),
 \]

which exists as a stochastic integral by the continuity of $f(t,s)$ in $s$.
The (standard) integration by parts formula yields for fixed $t$, thanks to the absolute continuity of $f$ in $s$ (Definition~\ref{d.1}(\ref{derivative_second_argument})),

\begin{equation}\label{e.integration_by_parts.1}
  \int_{\tau \lor -n }^t f(t,s) \ X( \dd s ) = f(t,t) X(t)-  f(t,\tau \lor -n ) X(\tau \lor -n ) - \int_{\tau \lor -n }^t X(s) \ds f(t,s) \ \dd s.
\end{equation}

In the case $\tau > -\infty$ we have $\tau \lor -n= \tau$ for $n$ sufficiently large, which proves~\eqref{integration_by_parts_2} in this case. If instead $\tau = -\infty$, we have
\[
\begin{aligned}
    |f(t,-n) X(-n)| 
  &=    |f(t,-n)\varphi(n)|\cdot \left|\frac{X(-n)}{\varphi(n)}\right|
  \\
  & \le  \sup_{s \in (-\infty,t]} \frac{|X(s)|}{\varphi(|s|)^{}}    |f(t,-n)\varphi(n)|.
\end{aligned}
\]
Since according to Lemma~\ref{e.exp.sup.finite} the first factor on the above right-hand side is bounded in $\mathscr L^p(\mathbb P)$ and therefore also bounded $\mathbb P$-a.s. and the second factor is deterministic and by \eqref{e.decay_of_f} tends to $0$ as $n \to \infty$, we deduce that 

\begin{equation}\label{e.convergence}
|f(t,-n) X(-n)|  \to 0 
\end{equation}

$\mathbb P$-a.s. and in $\mathscr L^p(\mathbb P)$ as $n \to \infty$. Moreover, since $s\mapsto\ds f(t,s)\varphi(|s|)\in\mathscr L^1(\lambda)$, cf. Remark~\ref{r.integral}, we obtain by Lemma~\ref{e.exp.sup.finite} that
%\begin{equation}\label{e.estimate_final}
\[
\int_{-\infty }^t \left|X(s) \ds f(t,s) \right| \ \dd s
\le 
\sup_{u\in(-\infty,t]}\left|\frac{X(u)}{\varphi(|u|)}\right|\int_{-\infty}^t \left| \ds f(t,s)\varphi(|s|)\right| \ \dd s<\infty
\]
$\mathbb P$-a.s. and in $\mathscr L^p(\mathbb P)$. 
Taking the limit as $n \to \infty$ in (\ref{e.integration_by_parts.1}) and using the dominated convergence theorem thus proves the assertion.

2. In view of Definition~\ref{d.1}(\ref{continuity_of_f}) and the assumption that $X$ is c\`adl\`ag we only have to show that the third term on the right-hand side of~\eqref{integration_by_parts_2} is continuous. 
For this purpose we set $\gamma := \frac{2\epsilon+1}{1+\epsilon} >1$, with $\epsilon>0$ as in Definition~\ref{d.1}(\ref{derivative_second_argument}), and define 
\[
\Upsilon(t,s):=X(s)\ds f(t,s)\left(1\lor |s|^{\gamma}\right)
\]
for any $s,t\in\R$. By means of \eqref{e.sup_ds} and Lemma~\ref{e.exp.sup.finite} we then obtain
\[
\begin{aligned}
 \sup_{t \in [0,T]} & 
 \int_{\tau}^t  \left|\Upsilon(t,s)\right|^{1+\epsilon} \ \frac{1}{1\lor|s|^\gamma}  \dd s \\
  &= \sup_{t \in [0,T]} \int_{\tau}^t  \left|X(s)   \ds f(t,s) \right|^{1+\epsilon} (1\lor|s|^{2\epsilon}) \  \dd s 
  \\
  &\le \sup_{t \in [0,T]} \int_{\tau}^t  \left|\varphi(|s|)   \ds f(t,s) \right|^{1+\epsilon} (1\lor|s|^{2\epsilon}) \  \dd s\cdot\sup_{s\in\R}\left(\frac{|X(s)|}{\varphi(|s|)}\right)^{1+\epsilon} <\infty
 \end{aligned}
\]
$\mathbb P$-almost surely. Consequently, we can use the de la Vall\'ee-Poussin theorem to deduce that $(\mathds1_{[\tau,t]}(\cdot)\Upsilon(t,\cdot))_{t \in [0,T]}$ is uniformly integrable with respect to the finite measure 
\[
\frac{1}{1\lor|s|^\gamma} \dd s.
\]
Now let $t \in [0,T]$ and choose an arbitrary sequence $(t_n)_{n \in \N}$ such that $t_n \to t$ as $n \to \infty$. The convergence $\Upsilon(t_n,s) \to \Upsilon(t,s)$ for $\lambda$-a.e. $s\in (-\infty,\tau]$, cf. Definition~\ref{d.1}\eqref{derivative_second_argument}(a), 
together with the uniform integrability of $(\mathds1_{[\tau,t]}(\cdot)\Upsilon(t,\cdot))_{t \in [0,T]}$ results in

\[
\begin{aligned}
 \lim_{n \to \infty}  \int_{\tau}^{t_n} X(s) \ds f(t_n,s) \ \dd s &=\lim_{n \to \infty}  \int_{\tau}^T \mathds1_{[\tau,t_n]}(s)\Upsilon(t_n,s)  \  \frac{1}{1\lor|s|^\gamma}  \dd s  
 \\
 &=  \int_{\tau}^t X(s) \ds f(t,s) \ \dd s.
 \end{aligned}
\]

This implies that the mapping $t \mapsto \int_{\tau}^t X(s) \ds f(t,s) \ \dd s$ is continuous.

3. Due to the c\`adl\`ag paths of $M$ (cf. 2.) the process $M$ is separable and thus $\sup_{t \in [0,T]}|M(t)|^p$ is measurable. 
By means of the integration by parts formula (\ref{integration_by_parts_2}), Minkowski's inequality, and Doob's inequality we obtain
\begin{eqnarray*}\label{e_supestimate}
 \left\| \sup_{t\in[0,T]}|M(t)| \,\right\|_p \nonumber  
&\le & \frac{p}{p-1} \sup_{t\in[0,T]}|f(t,t)|  \|X(T)\|_p
+ \sup_{t\in[0,T]} \|X(\tau)f(t,\tau)\|_p \nonumber \\ && +   \sup_{t\in[0,T]} \int_{\tau}^t \varphi(|s|) \left|\ds f(t,s)\right| \ \dd s
\cdot\left\|\sup_{s \in \mathbb{R}} \frac{|X(s)|}{\varphi(|s|)^{}} \right\|_p.
\end{eqnarray*}
The assertion is now a consequence of Definition~\ref{d.1}(\ref{continuity_of_f}), (\ref{e.convergence}), Remark~\ref{r.integral}, and Lemma~\ref{e.exp.sup.finite}.
\hfill$\square$

\section{The L\'evy-driven case}\label{s.cLp}

In this section we particularise the main results to the case when the driving martingale $X$ is a L\'evy process with 
some focus on fractional L\'evy processes as introduced by \cite{Mar06}.   

\subsection{Set-up}

Regarding some background on the theory of L\'evy processes we refer e.g. to
\cite{Ber96} and \cite{Kyp06}. Stochastic analysis with respect to L\'evy
processes is treated in \cite{App09}.

Let $(\gamma, \sigma, \nu)$ be a triplet consisting of constants $\gamma\in \R$ and $\sigma\ge0$ as well as a measure $\nu$ on $\R_0:=\R\setminus\{0\}$ that satisfy 
\begin{equation}\label{e.L-k}
\int_{\R_0} (x^2\land x) \ \nu(\dd x) <\infty
\end{equation}
and
\begin{equation}\label{drift_gamma}
  \gamma=-\int_{\R\setminus[-1,1]}x \ \nu(\dd x).
\end{equation}
Observe that $(\gamma, \sigma, \nu)$ is a so-called characteristic triplet  that determines the distribution of a L\'evy process on a complete probability space $(\Omega,\mathscr F,\mathbb P)$. Hence, let $\hat{L}:=(\hat{L}(t))_{t\geq 0}$ be a L\'evy process on  $(\Omega,\mathscr F,\mathbb P)$ with characteristic triplet $(\gamma, \sigma, \nu)$ and c\`adl\`ag paths, whose jump measure we denote by $N(\dd x, \dd s)$. Furthermore, let

%\begin{equation}\label{LK_representation}
\[
\begin{aligned}
  \Psi(u) := \ln\left(\mathbb E\left(e^{iu\hat{L}(1)}\right)\right)&=i\gamma u -\frac{\sigma^2 u^2}{2}+ \int_{\R_0} \left(e^{iux}-1- iux \mathds 1_{\left\{|x|\leq 1\right\}}\right) \  \nu(\dd x)\\
  &=-\frac{\sigma^2 u^2}{2}+ \int_{\R_0} \left(e^{iux}-1- iux \right) \  \nu(\dd x)
\end{aligned} 
\]

be the characteristic exponent of $\hat{L}$ which is given by the L\'evy-Khintchine formula. Note that (\ref{e.L-k}) implies that $\hat{L}(t)\in\mathscr L^1(\mathbb P)$ and
(\ref{drift_gamma}) is equivalent to $\E(\hat{L}(1))=0$. The latter assertion holds, since
\[
\E\left(\hat{L}(1)\right)=\Psi'(0+)=\int_\R \left.\frac{\dd}{\dd u}(e^{iux}- iux)\right|_{u=0} \  \nu(\dd x)=0,
\]
where $\Psi'$ denotes the derivative of $\Psi$. The process $\hat{L}(t)$ can be represented as

\[
\hat{L}(t) =  \sigma W(t)+ \int_0^t \int_{\R_0} x \ \tilde{N}(\dd x, \dd s), 
\]

where $W$ is a standard Brownian motion,  $\sigma$ is  the standard deviation of the Gaussian component of $\hat{L}(1)$, $\nu$ is the L\'evy measure and $ \tilde{N}(\dd x, \dd s) =  N(\dd x, \dd s)- \nu(\dd x) \dd s$ is the compensated jump measure of the L\'evy process $\hat{L}$.

By using two independent copies of $\hat L$ we consider a two-sided L\'evy process $L$ constructed in the spirit of (\ref{e.two-sided}). Having the process $L$ at hand we consider the \emph{L\'evy-driven Volterra process}

\begin{equation}\label{e.convoluted_levy_process}
  M(t) = \int_{-\infty}^t f(t,s) \ L(\dd s),
\end{equation}

which is a special case of \eqref{d.2}.
In the literature, processes as in \eqref{e.convoluted_levy_process} are occasionally also referred to as \emph{filtered L\'evy processes} (see e.g. \cite{DS_ac}) or \emph{convoluted L\'evy processes} (cf. \cite{BM_sccl}). However, as e.g. in \cite{BBV_mesp} we think that \emph{L\'evy-driven Volterra processes} is the most apposite name for such processes. In the special case that $f(t,s)=g(t-s)-g(-s)$ for some function $g$ such processes are also called \emph{moving average processes}.

The prime example of L\'evy-driven Volterra processes are \emph{fractional L\'evy processes}, where the integration kernel is given by
%\begin{equation}\label{e.fractional_kernel_definition}
\[
 f_d(t,s)=\frac{1}{\Gamma(d+1)}\left((t-s)^d_+-(-s)^d_+\right).
\]
Fractional L\'evy processes exist e.g. for parameters $d\in(0,1/2)$, when the driving centred L\'evy process $L$ is square integrable. In this case they have (up to some constant)
the same second order structure as a fractional Brownian motion with Hurst parameter $H=1/2+d$, but fractional L\'evy processes fail to be self-similar except in the fractional Brownian motion case, see \cite{Mar06}. 
The motivation for the name fractional L\'evy process is that it generalises the Mandelbrot-Van Ness representation of a fractional Brownian motion as an integral 
of the same kernel with respect to Brownian motion.

\subsection{Results on L\'evy-driven Volterra processes}

In order to make the result of Theorem \ref{t.0} applicable to the L\'evy-driven case, we need to control the $p$-th moment of $L$ as time approaches infinity.

\begin{lemma}\label{l.time_dependence_pth_moment}
 Let $p \ge 2$, $t \ge 1$ and $L(1) \in \L^p(\mathbb P)$. Then there exist a constant $C_p$ only depending on $p$ such that
 
 \[
  \E \left( |L(t)|^p\right) \le C_p t^{\frac{p}{2}} \E \left( |L(1)|^p\right).
 \]

\end{lemma}

\begin{pf}
 Note that for $l \ge 1$ and every L\'evy process $\tilde L$ with $\tilde L(1) \in \L^l(\mathbb P)$ and $n \in \N$ as well as $s \in [0,1)$ such that $t=n+s$ we infer by using Minkowski's inequality and the stationary increments of $\tilde L$
 
 \[
 \begin{aligned}
  \left\| \tilde L(t) \right\|_l &= \left\| \tilde L(s)+ \left(\tilde L(s+1)- \tilde L(s) \right)+ \ldots + \left(\tilde L(n+s)- \tilde L((n-1)+s) \right)\right\|_l \\
  &\le (t+1) \sup_{u \in [0,1]} \left\| \tilde L(u) \right\|_l.
 \end{aligned}  
 \]

 Since $t \ge 1$ implies $(t+1)^l \le (2t)^l$ the above leads to
 
 \begin{equation}\label{e.tilde_L_estimate}
  \E \left( \left|\tilde L(t)\right|^l\right)= \left\| \tilde L(t) \right\|_l^l  \le 2^l t^l \sup_{u \in [0,1]} \left\| \tilde L(u) \right\|_l^l.
 \end{equation}
 
 We now apply the Burkholder-Davis-Gundy inequality to the L\'evy process $L$ and deduce that there exists a constant $c_{p,1} > 0$ such that
 
 \begin{equation}\label{e.burkholder-davis-gundy}
    \E \left( \left| L(t)\right|^p\right) \le \E \left( \sup_{u \in [0,t]}\left| L(u)\right|^p\right) \le c_{p,1} \E \left( [L,L]_t^{\frac{p}{2}}\right).
 \end{equation}

Thus, choosing $\tilde L(t) := [L,L]_t$ and $l= \frac{p}{2}$, \eqref{e.tilde_L_estimate} and \eqref{e.burkholder-davis-gundy} result in

 \begin{equation}\label{e.1}
    \E \left( \left| L(t)\right|^p\right) \le c_{p,1} \E \left( [L,L]_t^{\frac{p}{2}}\right)\le c_{p,1} 2^{\frac{p}{2}}  \sup_{u \in [0,1]}\E \left( [L,L]_u^{\frac{p}{2}}\right) t^{\frac{p}{2}}= c_{p,1} 2^{\frac{p}{2}}  \E \left( [L,L]_1^{\frac{p}{2}}\right) t^{\frac{p}{2}}.
 \end{equation}

 We proceed with another application of the Burkholder-Davis-Gundy inequality and Doob's maximal inequality which lead to
 
 \[
  \E \left( [L,L]_1^{\frac{p}{2}}\right) \le c_{p,2} \E \left( \sup_{u \in [0,1]} |L(u)|^p \right) \le c_{p,2} \left(\frac{p}{p-1}\right)^p  \E \left(  |L(1)|^p \right)
 \]
 
 for some $c_{p,2}>0$. Plugging this into \eqref{e.1} results in 
 
 \[
  \E \left( \left| L(t)\right|^p\right) \le c_{p,1}c_{p,2} 2^{\frac{p}{2}}  \left(\frac{p}{p-1}\right)^p  \E \left( |L(1)|^p\right) t^{\frac{p}{2}}.
 \]
  \end{pf}

Theorem \ref{t.0} in the L\'evy-driven case roughly states that nice path and integrability properties of $L$ are carried over to the process $M$ for suitable kernel functions. The precise formulation reads as follows.

\begin{theorem}\label{t.1_levy}
Define
$\varphi_q(t)=|t|^q\vee 1$.
 Let $p \ge 2$ and suppose that $L$ is a centred L\'evy process such that $L(1)\in\mathscr L^p(\mathbb P)$. If $f\in\mathcal K(\varphi_q, -\infty)$ for some $q>1/2$, then 
$$
\tilde M(t)=\int_{-\infty}^t f(t,s)L(ds),\;0\leq t \leq T,
$$
exists as Wiener integral and a c\`adl\`ag modification of $\tilde M$ is given by 
$$
M(t) := f(t,t) L(t) - \int_{-\infty}^t L(s) \ds f(t,s) \ \dd s.
$$
This modification satisfies $\Delta M(t)=f(t,t)\Delta L(t)$ and the following maximal inequality:
There is a constant $C_{p,q}$ depending only on $(p,q)$ such that
\begin{eqnarray*}
\left\| \sup_{t \in [0,T]} |M(t)|\,\right\|_p\leq C_{p,q} \|L(1)\|_p\left( (T\vee 1)^{1/2} \sup_{t \in [0,T]} |f(t,t)|+\sup_{t \in [0,T]} \left(\int_{-\infty}^t  \left|\ds f(t,s)\right| (|s|^q\vee 1 )ds\right)\right).
\end{eqnarray*}
 \end{theorem}

{\bf Proof.} 
By the continuity of $f$ in the $s$-variable and its behaviour at $-\infty$ given by (\ref{e.decay_of_f}) with $\varphi=\varphi_q$ it is easy to check that
$f(t,\cdot)\in \mathscr L^2(\mathbb{R})$ for every $t$. Hence, $\tilde M(t)$ exists as a Wiener integral and
$$
\tilde M(t)=\lim_{n\rightarrow \infty} \int_{-n}^t f(t,s)L(ds)
$$
in $\mathscr L^2(\mathbb P)$. Using the stationary increments of $L$, Lemma~\ref{l.time_dependence_pth_moment} and the fact that $q>1/2$ we deduce

\begin{equation}\label{e.hilfc1}
 \sum_{n=0}^\infty \frac{\left\| L( 2^{n+1})- L(2^{n})\right\|_p}{\varphi_q(2^{n})} = \sum_{n=0}^\infty \frac{\left\| L( 2^{n})\right\|_p}{\varphi_q(2^{n})}  
\le C_p^{\frac{1}{p}} \E(|L(1)|^p)^{1/p} \sum_{n=0}^\infty (2^{n})^{1/2-q} < \infty.
\end{equation}
Hence condition \eqref{e.int_X_phi_finite} is satisfied and so Theorem \ref{t.0} applies with $L$ in place of $X$. 
Plugging (\ref{e.hilfc1}) into the right-hand side of (\ref{t.1}) and changing from $\|L(T)\|_p$ to $\|L(1)\|_p$ in (\ref{t.1})  via Lemma 
\ref{l.time_dependence_pth_moment} (if $T>1$) yields
the maximal inequality.
\hfill$\square$

\begin{example}\label{ex_fLp}
 We now come back to the case of a fractional L\'evy process and consider the kernel 
$$
f_d(t,s)=\frac{1}{\Gamma(d+1)}\left((t-s)^d_+-(-s)^d_+\right)
$$
for $d>0$. We first show that $f_d \in \mathcal K(\varphi_q, -\infty)$, if $d+q<1$. Indeed, by the mean value theorem we have for $t\in [0,T]$ and $s < 0$ 
\begin{eqnarray*}
 \Gamma(d+1) |f_d(t,s)|&\leq&  d t |s|^{d-1}, \\
 \Gamma(d+1)\left|\ds f_d(t,s)\right|&\leq& d(1-d) t |s|^{d-2}.
\end{eqnarray*}
The first inequality shows that (\ref{e.decay_of_f}) is satisfied with $\varphi=\varphi_q$, while the second one implies 
(\ref{e.sup_ds}) for $\epsilon<(1/(d+q)-1)\wedge (d/(1-d)).$
Hence, we observe that Theorem \ref{t.0} is applicable to fractional L\'evy processes $M_d$ for $0<d<1/2$, if the driving centred L\'evy process has a finite second moment. Continuity of
the fractional L\'evy process follows from $f_d(t,t)=0$, but is well-known (see e.g. \cite{Mar06}). By the substitution $s=vt$ we obtain for $q>1/2$ such that $d+q<1$
 \begin{eqnarray*}
  && \Gamma(d+1) \int_{-\infty}^t  \left|\ds f_d(t,s)\right| (|s|^q\vee 1 )ds\\ &\leq& 
t^{d+q}\int_{-\infty}^1 \left((1-v)_+^{d-1}-(-v)_+^{d-1}\right) |v|^q dv + t^{d} \int_{-\infty}^1 \left((1-v)_+^{d-1}-(-v)_+^{d-1}\right) dv.
 \end{eqnarray*}
Hence, the maximal inequality in Theorem~\ref{t.1_levy} can be simplified as follows: If  $q>1/2$ and  $d+q<1$, there is a constant $C_{p,q,d}$ independent of the driving L\'evy process $L$ such that for every $T\geq 1$
$$
\left\| \sup_{t \in [0,T]} |M_d(t)|\,\right\|_p\leq C_{p,q,d} \|L(1)\|_p \; T^{d+q}.
$$
\end{example}

We finally consider the case when $L$ is a symmetric stable process with index of stability $\alpha\in (1,2)$. In this case the L\'evy-driven Volterra processes generalise fractional $\alpha$-stable motions,
for which we refer to~\cite{ST_sngrp} and \cite{EW_uaflp}. The following variant of Theorem \ref{t.1_levy} holds true.
 
\begin{theorem}\label{t.2_levy}
Suppose that $L$ is a symmetric $\alpha$-stable L\'evy process with $\alpha\in (1,2)$.  If $f\in\mathcal K(\varphi_q, -\infty)$ for some $q>1/\alpha$, then 
$$
\tilde M(t)=\int_{-\infty}^t f(t,s)L(ds),\;0\leq t \leq T,
$$
exists as an $L$-integral in the sense of \cite{RR89} and a c\`adl\`ag modification of $\tilde M$ is given by 
$$
M(t) := f(t,t) L(t) - \int_{-\infty}^t L(s) \ds f(t,s) \ \dd s.
$$
This modification satisfies $\Delta M(t)=f(t,t)\Delta L(t)$ and the following maximal inequality:
For every $p\in (1,\alpha)$ there is a constant $C_{p,q}$ depending only on $(p,q)$ such that
\begin{eqnarray*}
\left\| \sup_{t \in [0,T]} |M(t)|\,\right\|_p\leq C_{p,q} \|L(1)\|_p\left( T^{1/\alpha} \sup_{t \in [0,T]} |f(t,t)|+\sup_{t \in [0,T]} \left(\int_{-\infty}^t  \left|\ds f(t,s)\right| (|s|^q\vee 1 )ds\right)\right).
\end{eqnarray*}
 \end{theorem}
\begin{pf}
 We first show existence of the $L$-integral. As $f \in K(\varphi_q, -\infty)$, there is a constant $C(t)>0$ depending on $t$ such that
$$
|f(t,s)|\leq C(t)(|s|^{-q}\wedge 1)
$$
for every $s\in \mathbb{R}$.
Noting that the L\'evy measure of a symmetric $\alpha$-stable process is given by $\nu(dx)=A |x|^{-1-\alpha} dx$ for some constant $A>0$, we get for every $t\in [0,T]$, $1<p<\alpha$, and $2-\alpha<\gamma<2-1/q$
\begin{eqnarray*}
 && \int_{\mathbb{R}} \int_{\mathbb{R}_0} \left( |f(t,s)x|^2{\mathds{1}}_{\{|f(t,s)x|\leq 1\}} +  |f(t,s)x|^p{\mathds{1}}_{\{|f(t,s)x|> 1\}}\right) \nu(dx) ds
\\ &\leq &  A \int_{\mathbb{R}} |f(t,s)|^{2-\gamma}ds \int_{\{|x|\geq 1\}} |x|^{1-\gamma-\alpha} dx + A \int_{\mathbb{R}} |f(t,s)|^{2} ds \int_{\{0<|x|<1\}} |x|^{1-\alpha} dx
\\ && + A\int_{\mathbb{R}} |f(t,s)|^p \left(\int_{\{|x|\geq C(t)^{-1} (|s|^q \vee 1)\}} |x|^{p-1-\alpha} dx \right) ds
\\ &\leq & \tilde C(t) \int_{\mathbb{R}}  \left( (|s|^{-q(2-\gamma)}\wedge 1) + (|s|^{-2q}\wedge 1) +  (|s|^{-q\alpha}\wedge 1)\right)ds <\infty,
\end{eqnarray*}
where the constant $\tilde C(t)$ depends on $\gamma$ and $t$. Taking the symmetry of $\nu$ into account, Theorem 3.3 in~\cite{RR89} implies that 
$\tilde M(t)$ exists as an $L$-integral and that
$$
\tilde M(t)=\lim_{n\rightarrow \infty} \int_{-n}^t f(t,s)L(ds)
$$
in $\mathscr L^p(\mathbb P)$ .

The remainder of the proof is analogous to that of Theorem \ref{t.1_levy}. Instead of Lemma~\ref{l.time_dependence_pth_moment} one can apply that 
$$
\|L(t)\|_p=t^{1/\alpha} \|L(1)\|_p
$$
for  all $t\geq 0$ by the self-similarity of the symmetric $\alpha$-stable process.
\end{pf}

\section*{Acknowledgements}

The authors are grateful to Alexander Lindner for his help to extend Lemma~\ref{l.time_dependence_pth_moment} from the case of even integers to the general case.
Financial support by the Deutsche Forschungsgemeinschaft under grant BE3933/4-1 is gratefully acknowledged.

% %\bibliographystyle{unsrt}
% \bibliographystyle{plain}  
% \bibliography{literatur}
\end{document}